\RequirePackage{fix-cm}

\documentclass[smallextended]{svjour3} 
\smartqed  

\usepackage{amssymb,amsmath,amsfonts,latexsym} 
\usepackage{hyperref}
   
\def\proofend{\hfill$\Box$} 

\DeclareMathOperator{\sgn}{sgn} 

\usepackage{mathtools}
\usepackage{multirow,array}
\usepackage{graphicx}
\usepackage{subfigure}
\usepackage{epstopdf}
\usepackage{cite}
\usepackage{float} 
\usepackage{color, xcolor}
\allowdisplaybreaks

\usepackage{calrsfs}
\DeclareMathAlphabet{\pazocal}{OMS}{zplm}{m}{n}

\newcommand{\be}{\begin{equation}}
\newcommand{\ee}{\end{equation}}
\newcommand{\RR}{\mathbb{R}}
\newcommand{\CC}{\mathbb{C}}
\newcommand{\NN}{\mathbb{N}}
\newcommand{\mL}{\mathcal{L}}
\newcommand{\Leb}{\pazocal{L}}
\newcommand{\mD}{\mathcal{D}}
\newcommand{\DC}{\mathbb{D}}
\newcommand{\DRL}{\pazocal{D}}
\newcommand{\mE}{\mathcal{E}}
\newcommand{\mF}{\mathcal{F}}
\newcommand{\mM}{\mathcal{M}}
\newcommand{\mP}{\mathcal{P}}

\journalname{Fract. Calc. Appl. Anal.} 

\begin{document}


\title{Mellin definition of the fractional Laplacian}

\titlerunning{Mellin definition of the fractional Laplacian}

\author{Gianni Pagnini$^{1,2}$ 
\and
Claudio Runfola$^{3}$ 
}

\authorrunning{G. Pagnini \and C. Runfola} 

\institute{Gianni Pagnini$^{1,2,*}$
\at
$^1$BCAM -- Basque Center for Applied Mathematics, 
Alameda de Mazarredo 14, E-48009 Bilbao, Basque Country -- Spain \\
$^2$Ikerbasque -- Basque Foundation for Science, 
Plaza Euskadi 5, E-48009 Bilbao, Basque Country -- Spain \\
\email{gpagnini@bcamath.org}
$^*$ corresponding author \\ 
ORCID ID: orcid.org/0000-0001-9917-4614 
 \and
Claudio Runfola$^{3}$
\at
$^3$Independent Researcher\\
\email{claudio.runfola@gmail.com} 
}

\date{Received: / Revised: / Accepted: ......}

\maketitle

\begin{abstract}
It is known that
at least ten equivalent definitions of the fractional Laplacian
exist in an unbounded domain.
Here we derive a further equivalent definition that is based on
the Mellin transform and it can be used when 
the fractional Laplacian is applied to radial functions.
The main finding is tested in the case of 
the space-fractional diffusion equation.
The one-dimensional case is also considered, 
such that the Mellin transform of the Riesz 
(namely the symmetric Riesz--Feller)
fractional derivative is established. 
This one-dimensional result corrects an
existing formula in literature. 
Further results for the Riesz fractional derivative are obtained
when it is applied to symmetric functions,
in particular its relation with the Caputo and the Riemann--Liouville
fractional derivatives.


\end{abstract} 


\section{Introduction} \label{sec:1}
\setcounter{section}{1} 
\setcounter{equation}{0} 

The fractional Laplacian \cite{bucur_valdinoci-2016,lischke_etal-jcp-2020}, 
or fractional Laplace operator, is a relevant example of
non-local pseudo-differential operator 
within the theory of fractional calculus.
In an unbounded domain,
there exist at least ten equivalent definitions 
of the fractional Laplacian \cite{kwasnicki-fcaa-2017},
while this equivalence does not hold in bounded domains,
see, e.g., 
\cite{servadei_etal-prse-2014,cusimano_etal-siamjna-2018,
duo_etal-dcdsb-2019}. 
Here we show that a further equivalent definition in 
unbounded domains can be 
provided on the basis of the Mellin transform
when the fractional Laplacian is applied to radial functions.
Actually, the Mellin transform has a remarkable role
in fractional calculus, see, e.g., 
\cite{gorenflo_etal-fcaa-2000,mainardi_etal-fcaa-2001,mainardi_etal-fcaa-2003,
luchko_etal-fcaa-2013,bardaro_etal-jfaa-2015},
that is additionally supported by the present results.

More in detail, this equivalence holds when the fractional Laplacian
is an operator on Lebesgue spaces $\Leb^p$ (with $p \in [1,\infty)$), 
on the space of continuous functions vanishing at infinity and 
on the space of bounded uniformly continuous functions
\cite{kwasnicki-fcaa-2017}.  
For the present aims, and to lighten the text,
we assume sufficiently well-behaving functions 
in $\RR^n$ with $n \in \{1,2,\dots\}$ such that 
the considered formulae have sense.
In particular, 
we consider radial functions in $\RR^n$ that are 
radial Schwartz--Bruhat functions, 
namely they and all their derivatives vanish very quickly for 
large value of the argument such that we can use a result previously derived
by N. Ormerod \cite{ormerod-jmaa-1979}. 
We denote the fractional Laplacian in $\RR^n$ 
by $\mL=-(-\Delta)^{\alpha/2}$ with $\alpha \in (0,2)$.
The classical Laplacian $\Delta$ is recovered
when $\alpha=2$. 
The reader interested in the technicalities concerning the 
fractional Laplacian is referred,
for example, to the book by C. Bucur and E. Valdinoci 
\cite{bucur_valdinoci-2016}, 
and for the technicalities concerning the Mellin transform
is referred, for example, 
to the book by O.I. Marichev \cite{marichev-1983} 
or to the book by R.B. Paris and D. Kaminski \cite{paris_kaminski-2001}. 

Before to start, we want to remind also that 
the fractional Laplacian is a peculiar nonlocal operator
beside its role in fractional calculus. 
In particular, it has been proved that
all functions are locally $s$-harmonic up to a small error
\cite{dipierro_etal-jems-2017} in the sense
of the fractional Laplacian of order $\alpha/2=s$. 
Actually,
this is a pure consequence of nonlocality of the fractional Laplacian
because the approximating functions that are 
$s$-harmonic are defined in the whole $\RR^n$ such that 
the approximation is possible with a careful choice of the values
outside the reference interval \cite{krylov-jfa-2019}. 
A similar property does not hold for any nonlocal operator 
\cite{dipierro_etal-jems-2017}. 
This result has been further investigated in literature
\cite{valdinoci-lnm-2018,dipierro_etal-jga-2019,krylov-jfa-2019}.

We begin by reminding the heat semi-group definition 
\cite{kwasnicki-fcaa-2017,stinga-2019,lischke_etal-jcp-2020} 
(or Bochner's definition \cite{bochner-1949})
\be
\mL u = \frac{-1}{\Gamma\left(-\frac{\alpha}{2}\right)}
\int_0^\infty (e^{t\Delta}u -u) \frac{dt}{t^{1+\alpha/2}} \,,
\quad \alpha \in (0,2) \,,
\label{eq:semigroup}
\ee
where $e^{t\Delta}$ is the propagator of the heat equation,
which means that 
$\rho(x,t)=e^{t\Delta}u(x):\RR^n \times [0,\infty) \to \RR$ solves
\be
\frac{\partial \rho}{\partial t} = \Delta \rho \,, 
\quad \rho(x,0)=u(x) \,.
\label{eq:heat}
\ee
Then, by applying the Fourier transform to \eqref{eq:semigroup}, 
we have another and equivalent definition of the 
fractional Laplacian. In particular, 
the fractional Laplacian emerges to be a Fourier multiplier 
with the celebreted symbol
\be
\mF(\mL u)(\kappa)=
- |\kappa|^\alpha \mF u(\kappa) \,, \quad
\kappa \in \RR^n \,.
\label{eq:fourier}
\ee
Formula \eqref{eq:fourier} can be obtained by noting
from \eqref{eq:heat} that 
$\mF \rho(\kappa) = e^{-|\kappa|^2 t} \mF u(\kappa)$,
where $|\kappa|=\sqrt{\kappa \cdot \kappa}$ and
$|\kappa|^\alpha=(|\kappa|^2)^{\alpha/2}$, 
and by remembering that 
it holds
$\displaystyle{\int_0^\infty(e^{-y}-1) \, y^{s-1} \, dy=\Gamma(s)}$
for $s \in \CC$ 
with $-1 < \Re(s) < 0$
\cite[Chapter XII]{whittaker_watson-1952}
from the analytical continuation of the Gamma function 
through the Cauchy--Saalsch{\"u}tz representation,
see also \cite[Chapter 3]{temme-1996} and
\cite[Table of Mellin transforms]{paris_kaminski-2001}.

The fractional Laplacian is also related 
with the so-called Riesz potential $I^\alpha$ that is defined as 
\be
I^\alpha u(x) = 
\frac{1}{\gamma_n(\alpha)}
\int_{\RR^n} u(x+z)|z|^{-n+\alpha} \, dz \,,
\quad \gamma_n(\alpha)=\frac{\Gamma\left(\frac{n-\alpha}{2}\right)}
{2^\alpha \pi^{n/2} \Gamma \left(\frac{\alpha}{2}\right)} \,,
\label{RieszI}
\ee
with $\alpha \ne n + 2m$, $m=0,1,2, \dots$, and converges
when $0 < \alpha < n$, 
that is with $\alpha \in (0,2)$ when $n \ge 2$ 
and $\alpha \in (0,1)$ when $n =1$,
see, e.g., \cite{samko-fcaa-1998,rubin-fcaa-2012}. 
In fact, another equivalent definition of the fractional Laplacian 
follows from the left inverse of the Riesz potential 
\eqref{RieszI}
\cite{samko-fcaa-1998,kwasnicki-fcaa-2017}.
Actually, integral \eqref{RieszI} was originally introduced
by M. Riesz in the late 1940's 
in a study on the Cauchy problem for the wave equation 
with the aim of writing the solution in closed form
\cite{riesz-1949}. As a matter of fact,  
the operator \eqref{RieszI} defines 
a family of convolution operators with fractional index. 

In this paper we prove a further 
equivalent definition
of the fractional Laplacian on the basis of the Mellin transform.

When $n=1$, the left inversion of the Riesz potential \eqref{RieszI}
is called Riesz fractional derivative and it is here denoted by 
$\mD_x^\alpha$ \cite{cai_etal-fcaa-2019}. 
Thus, in the one-dimensional case, 
the symbol of the fractional Laplacian \eqref{eq:fourier}
and that of the Riesz fractional derivative are identical
because only the modulus of the Fourier variable is involved.
The Riesz fractional derivative can be written by highlighting its relation
with the second derivative, i.e., the one-dimensional Laplacian, 
through the following equality involving their symbols 
$-|\kappa|^\alpha = - \left(|\kappa|^2\right)^{\alpha/2}$ 
\cite{mainardi_etal-fcaa-2001} 
and then, in analogy, the following notation can be formally adopted
$\displaystyle{
\mF^{-1} (\mF \mD_x^\alpha) = \mD_x^\alpha 
:=  - \left( - d^2/dx^2 \right)^{\alpha/2}}$
when $x \in \RR$ \cite{feller-1952}. 
Moreover, we report here that W. Feller generalised the Riesz potential 
\eqref{RieszI} when $n=1$ to any complex $\alpha$ with
$\Re(\alpha) > 0$
and included a further parameter $\theta$ related with 
operator asymmetry \cite{feller-1952}.
In particular, by following the approach due to M. Riesz, 
W. Feller introduced a family of one-dimensional pseudo-differential 
operators obtained 
by the inversion of linear combinations of left and 
right hand-sided Riemann--Liouville--Weyl operators 
\cite{feller-1952,gorenflo_etal-fcaa-1998}. 
Furthermore, 
motivated by the research started by S. Bochner \cite{bochner-1949}
on the generalization of standard diffusion 
to generalized diffusion equations for L\'evy stable densities,
W. Feller provided such generalized operators,
nowadays called Riesz--Feller fractional derivative,
for diffusion processes converging to L\'evy stable densities 
on the basis of the fact that
each stable density defines a semi-group of transformations \cite{feller-1952}. 
Finally, the symbol of the Riesz--Feller fractional derivative
is $-|\kappa|^\alpha e^{i(\sgn \kappa)\theta\pi/2}$ 
with $\alpha \in (0,2)$ and $|\theta|\le\min\{\alpha,2-\alpha\}$. 
Thus, the symmetric Riesz--Feller fractional derivative
is the Riesz fractional derivative when $\theta=0$,
that is when the Riesz--Feller derivative is a symmetric operator,
and it is consistent with the requirement of radial functions in $\RR^n$.
Actually, the Riesz fractional derivative, 
or more general two-sided fractional derivatives,
are still topics for recent research, see, e.g.,
\cite{bayin-jmp-2016,cai_etal-fcaa-2019,ortigueira-mmas-2021}.
The Riesz fractional derivative 
can be written explicitly in a regularised integral
form valid for $\alpha \in (0,2)$ 
\cite[see, formula (2.6) and footnote $3$]{mainardi_etal-fcaa-2001}
\be 
\mD_x^\alpha f(x)=
\frac{\Gamma(1+\alpha)}{\pi} \sin\frac{\pi\alpha}{2}
\int_0^\infty
\frac{f(x+z)-2f(x)+f(x-z)}{z^{1+\alpha}} \, dz \,.
\label{RieszD}
\ee

By exploiting the one-dimensional case of
the Mellin definition of the fractional Laplacian here derived,
we also establish the relations of 
the Riesz fractional derivative applied to a symmetric function
with the Caputo and with the Riemann--Liouville fractional derivatives.

In the following, 
we first present and discuss the main results in Section \ref{sec:2}
and later we provide the proofs in Section \ref{sec:3}.
In Section \ref{sec:4}, 
we apply (and check) our findings for solving the space-fractional
diffusion equation. 
Conclusions and future developments
are finally reported in Section \ref{sec:5}.


\section{Main results}
\label{sec:2}
\setcounter{section}{2} 
\setcounter{equation}{0} 

\begin{theorem}\label{Th1} 
Let $f$ be a radial function in $\RR^n$
and $\mM(\mL f)(s)$ its Mellin transform with $s \in \CC$.
A further definition of the fractional Laplacian 
$\mL=-(-\Delta)^{\alpha/2}$ with $\alpha \in (0,2)$
that is equivalent, inter alia, 
to definitions \eqref{eq:semigroup} and \eqref{eq:fourier} 
in the sense of Reference \cite{kwasnicki-fcaa-2017} is
\be
\mM(\mL f)(s)=
- 2^\alpha 
\frac{\Gamma\left(\frac{s}{2}\right) 
\Gamma\left[\frac{n-(s-\alpha)}{2}\right]}
{\Gamma\left(\frac{n-s}{2}\right) 
\Gamma\left[\frac{s-\alpha}{2}\right]}
\mM f(s-\alpha) \,,
\quad s \in \CC \,,
\label{eq:th1}
\ee
for $s \ne 0$ and $0 < \Re(s)<n$.
\end{theorem}

From Theorem \ref{Th1} with $n=1$, 
we have the following Corollary that
holds for the Riesz fractional derivative \eqref{RieszD}

\begin{corollary}\label{Cor1} 
When the Riesz frational derivative
$\mD_x^\alpha$ \eqref{RieszD} of order $\alpha \in (0,2)$ 
is applied to a symmetric function $f(x)=f(|x|)$, with $x \in \RR$,
it holds
\be
\mM(\mD_x^\alpha f)(s)=
- \frac{\Gamma(s)\cos\left(\frac{\pi}{2}s\right)}
{\Gamma(s-\alpha)\cos\left[\frac{\pi}{2}(s-\alpha)\right]} 
\mM f(s-\alpha) \,,
\quad s \in \CC \,,
\label{eq:corollary}
\ee
for $s \ne 0$ and $0 < \Re(s)<1$.
\end{corollary}

Formula \eqref{eq:corollary} corrects a similar 
functional equation previously derived in an attempt to derive 
a new generalization of an Hilfer-type fractional derivative
\cite[formula (53)]{khan_etal-rp-2021}.
There \cite{khan_etal-rp-2021}, the expression 
of the Mellin transform of the Riesz fractional derivative 
is retrieved by handling this last as a pseudo-differential 
operator of the Riemann--Liouville type \cite[formula (42)]{khan_etal-rp-2021}. 
However, that formula does not work when it is faced 
to the space-fractional diffusion equation. 
The success of formula \eqref{eq:th1} 
for solving the space-fractional diffusion equation, 
and then also of \eqref{eq:corollary} for the one-dimensional case, 
is indeed here reported.

Moreover, the study of the one-dimensional case
allows for deriving 
the following functional relations that are
of interest in the theory of fractional calculus.
Let $x\in\RR$ and 
let $\DC_{|x|}^\alpha$ and $\DRL_{|x|}^\alpha$  
be the fractional derivatives of order $\alpha > 0$ 
with respect to $|x|\in\RR^+$
in the Caputo and in the Riemann--Liouville sense, respectively,
then the following Theorems can be proven
by starting from results obtained from Corollary \ref{Cor1}.

\begin{theorem}
\label{Th:DRDC}
Let $\phi:\RR \to \RR$ be a symmetric function such that
$\phi(x)=\phi(|x|)$ and 
let $\mD_x^\alpha$ and $\DC_{|x|}^\alpha$ be the Riesz and the
Caputo fractional derivatives with respect to $x$ and $|x|$, respectively,
of order $0\le m-1 < \alpha < m$ with $m\in\NN$.
Provided that 
\be
i) \, \lim_{|x| \to 0} \phi^{(j)}(|x|) \, |x|^{s-\alpha+j} = 0 \,, \quad
ii) \, \lim_{|x| \to +\infty} \phi^{(j)}|x|^{s-\alpha+j} = 0 \,, 
\label{cond:DRDC}
\ee
with $j=0,\dots,m-1$ and proper $\Re(s)$, it holds 
\be
\mD_x^\alpha \phi(|x|)=
\int_0^\infty \left\{\DC_{\xi}^\alpha \phi(\xi)\right\} 
\mE(|x|/\xi) \, \frac{d\xi}{\xi} \,,
\label{eq:DRDC}
\ee
where, for $\xi \in \RR^+$, kernel $\mE(\xi)$ is defined as
\[
\mE(\xi)= \frac{2}{\pi} 
\frac{\sin(\pi\alpha/2)}{1-\xi^2} - 
\cos\left(\frac{\pi\alpha}{2}\right) \, \xi \delta(\xi-1) \,.
\]
\end{theorem}

\begin{corollary}
When $\alpha=1$, it holds $\DC_{\xi}^1=d/d\xi$ with $\xi\in\RR^+$. 
Thus,
by using the Leibniz integral rule and the symmetry of $\phi$ 
such that $d\phi/dx=\sgn(x) \, d\phi/d|x|$ with $x\in\RR$, 
from Theorem \ref{Th:DRDC} we have
\be
\mD_x^1 \phi(|x|) = 
-\frac{1}{\pi}\frac{d}{dx} 
\int_{-\infty}^{+\infty} \frac{\phi(|y|)}{x-y} \, dy \,.
\label{eq:hilbert}
\ee
Formula \eqref{eq:hilbert} was previosuly derived
by W. Feller in the special case of the symmetric Cauchy density
by looking for space-fractional diffusion equations solved by the
stable densities introduced by P. L\'evy \cite{feller-1952}.
\end{corollary}
\begin{remark}
Formula \eqref{eq:hilbert} holds for symmetric functions only.
This symmetry constraint  
is implicit from the symmetric diffusion process studied
in the original derivation by W. Feller \cite{feller-1952}, 
but it is indeed omitted 
when the operator $\mD_x^1$ is reported in overviews, see,
e.g., \cite{mainardi_etal-fcaa-2001,bayin-jmp-2016}.
On the other side, formula \eqref{eq:hilbert} is an important
complement that is missing in recent studies on 
the Riesz fractional derivative or on more general
two-sided fractional derivatives \cite{cai_etal-fcaa-2019,ortigueira-mmas-2021}.
\end{remark}
\begin{corollary} 
When $\alpha=2$, by using the properties of the delta-function, 
from Theorem \ref{Th:DRDC} it holds
\be
\mD_x^2 \phi(|x|)=
\DC_{|x|}^2 \phi(|x|) \,.
\label{eq:DRDC2}
\ee
\end{corollary}

\begin{theorem}
\label{Th:DRDRL}
Let $\phi:\RR \to \RR$ be a symmetric function such that
$\phi(x)=\phi(|x|)$ and 
let $\mD_x^\alpha$ and $\DRL_{|x|}^\alpha$ be the Riesz and the
Riemann--Liouville fractional derivatives with respect to $x$ and $|x|$, 
respectively, of order $0\le m-1 < \alpha < m$ with $m\in\NN$.
Provided that conditions $i)$ and $ii)$ in \eqref{cond:DRDC} are met
together with the further conditions $\phi^{(j)}(0^+)=0$, 
for any $j=0,\dots,m-1$, then it holds 
\be
\mD_x^\alpha \phi(|x|)=
\int_0^\infty \left\{\DRL_{\xi}^\alpha \phi(\xi)\right\} 
\mE(|x|/\xi) \, \frac{d\xi}{\xi} \,,
\label{eq:DRDRL}
\ee
where kernel $\mE$ is defined as in Theorem \ref{Th:DRDC}.
\end{theorem}

\begin{corollary}
A further Theorem similar to Theorems \ref{Th:DRDC} and \ref{Th:DRDRL}
can be proved for the Miller--Ross fractional derivative.
Actually, a formula analogous to formuale \eqref{eq:DRDC} and \eqref{eq:DRDRL} 
with the same kernel $\mE$ \eqref{eq:E} can be derived after providing
the conditions required for the Mellin transform of the
Miller--Ross fractional derivative to be equal to that 
of both the Caputo and the Riemann--Liouville fractional
derivatives \cite[Section 2.10.5]{podlubny-1999}.
\end{corollary}

\section{Proofs and further results} 
\label{sec:3}
\setcounter{section}{3} 
\setcounter{equation}{0} 

\subsection{Proof of Theorem \ref{Th1}}
Before proving our Theorem \ref{Th1} we need to remind 
a previous result relating Fourier and Mellin transforms
of radial functions in $\RR^n$. 
This result was originally proved by N. Ormerod 
\cite[Theorem 1]{ormerod-jmaa-1979}.
Actually, N. Ormerod derived this as an intermediate step for providing 
an alternative proof of an existing result linking
Fourier transform of radial function in $\RR^n$ with 
appropriate Bessel transform that was previously derivided by using
different techniques \cite{bochner_chadrasekharan-1949,stein_weiss-1972}. 
However, 
it is this intermediate step in terms of Mellin transform 
provided by N. Ormerod that is useful for our proof.

Let the Fourier transform of a suitable function 
$g:\RR^n \to \RR$ be defined as
\be
\mF g(\kappa)=
\int_{\RR^n} e^{i \kappa \cdot x} \, g(x) \, dx \,,
\quad \kappa \in \RR^n \,,
\label{eq:fouriert}
\ee
then we have formula \eqref{eq:fourier}. 
Moreover, let the Mellin transform of another suitable function 
$\varphi:\RR_0^+ \to \RR$ be defined as 
\be
\mM \varphi(s)=\int_0^\infty \varphi(r) \, r^{s-1} \, dr \,,
\quad s \in \CC \,.
\label{eq:mellint}
\ee  
Then, let $f$ be a radial function in $\RR^n$,
i.e., $f(x)=f(r)$ where $r=|x| \in \RR_0^+$, 
in our notation the Theorem by Ormerod reads:

\begin{theorem}[Ormerod \cite{ormerod-jmaa-1979}]
\label{Th2} 
Function $\mM f(s)$ has an analytic continuation valid for
all $s \ne 0$ and satisfies the functional equation
\be
\mM f(s)=
\pi^{n/2-s}  
\frac{\Gamma\left(\frac{s}{2}\right)} 
{\Gamma\left(\frac{n-s}{2}\right)} 
\mM (\mF' f) (n-s) \,,
\label{eq:th2}
\ee
where for a suitable function $g$ the operator $\mF'$
is defined as $\mF'g(\kappa)=\mF g(2\pi\kappa)$.
\end{theorem}
More specifically, Theorem \ref{Th2} allows for the identity
\cite[equation (3)]{ormerod-jmaa-1979}
\be
\int_0^\infty f(r) \, r^{s-1} \, dr=
\pi^{n/2-s}  
\frac{\Gamma\left(\frac{s}{2}\right)} 
{\Gamma\left(\frac{n-s}{2}\right)} 
\int_0^\infty \{\mF' f(r)\} \, r^{n-s} \, dr \,,
\label{eq:identity}
\ee
for $s \ne 0$ and $0 < \Re(s) < n$.
Now, we can prove Theorem \ref{Th1}.  

\proof 
By using formula \eqref{eq:th2} from Theorem \ref{Th2}, 
we can write the functional equation
\be
\mM(\mL f)(s)=
\pi^{n/2-s}  
\frac{\Gamma\left(\frac{s}{2}\right)} 
{\Gamma\left(\frac{n-s}{2}\right)} 
\mM \{\mF'(\mL f)\} (n-s) \,.
\label{eq:proof1}
\ee
By adopting the operator $\mF'$, 
formula \eqref{eq:fourier} for a radial function $f$ reads
$\mF'(\mL f)(\kappa)=
- (2\pi |\kappa|)^\alpha \mF' f(|\kappa|)$.
Hence, the term $\mM \{\mF'(\mL f)\} (n-s)$ in \eqref{eq:proof1}
results to be
\begin{eqnarray}
\mM \{\mF'(\mL f)\}(n-s) 
&=& - (2\pi)^\alpha 
\int_0^\infty |\kappa|^\alpha \mF' f(|\kappa|) |\kappa|^{n-s-1} \, d|\kappa|
\nonumber \\
&=& - (2\pi)^\alpha \mM \{\mF'f\}(n-(s-\alpha)) 
\nonumber \\
&=& - 2^\alpha \pi^{-n/2+s}  
\frac{\Gamma\left[\frac{n-(s-\alpha)}{2}\right]} 
{\Gamma\left(\frac{s-\alpha}{2}\right)} 
\mM f(s-\alpha) \,,
\label{eq:proof2}
\end{eqnarray}
where Theorem \ref{Th2} has been used 
for calculating the term $\mM \{\mF'f\}(n-(s-\alpha))$
in the last step after replacing $s$ with $s-\alpha$.
Finally, by plugging the last line of \eqref{eq:proof2}
into \eqref{eq:proof1}, 
and by reminding the conditions for the validity of
\eqref{eq:th2} and identity \eqref{eq:identity}, we have formula \eqref{eq:th1}.
This completes the proof.
\proofend 

\subsection{Proof of Corollary \ref{Cor1}}
Before to proove Corollary \ref{Cor1},
we remind the following properties of the Gamma function
\be
\Gamma(2z)=\frac{2^{2z-1/2}}{\sqrt{2\pi}}
\Gamma(z)\Gamma\left(z+\frac{1}{2}\right) \,,
\label{eq:gammaproperty1}
\ee
\be
\Gamma(z)\Gamma(1-z)=\frac{\pi}{\sin \pi z} \,,
\label{eq:gammaproperty2}
\ee
and the trigonometric property
\be
\sin\left[\frac{\pi}{2}(1-\phi)\right] =
\cos\left[\frac{\pi}{2}\phi\right] \,.
\label{eq:trig}
\ee
By using formulae \eqref{eq:gammaproperty1}, 
\eqref{eq:gammaproperty2} and \eqref{eq:trig} we
can derive formula \eqref{eq:corollary} from formula \eqref{eq:th1}
when $n=1$.

\proof
By applying \eqref{eq:gammaproperty1} 
to the terms $\Gamma(s/2)$ and $\Gamma[(s-\alpha)/2]$
in formula \eqref{eq:th1} this last reads 
\be
\mM(\mL f)(s)=
- \frac{\Gamma(s)}{\Gamma(s-\alpha)} 
\frac{\Gamma\left(\frac{s-\alpha+1}{2}\right) 
\Gamma\left[\frac{n-(s-\alpha)}{2}\right]}
{\Gamma\left(\frac{s+1}{2}\right) 
\Gamma\left[\frac{n-s}{2}\right]}
\mM f(s-\alpha) \,,
\label{eq:formula}
\ee
for $s \in \CC$ with $s \ne 0$ and $0 < \Re(s) < n$.
By setting $n=1$ and using properties 
\eqref{eq:gammaproperty2} and \eqref{eq:trig},
for the numerator of \eqref{eq:formula} it holds
\begin{eqnarray}
\Gamma\left[\frac{1-(s-\alpha)}{2}\right]
\Gamma\left[\frac{1+(s-\alpha)}{2}\right]
&=&
\Gamma\left[\frac{1-(s-\alpha)}{2}\right]
\Gamma\left[1-\frac{1-(s-\alpha)}{2}\right] \nonumber \\
&=& \frac{\pi}{\sin\left[\pi(1-(s-\alpha))/2\right]} \nonumber \\
&=& \frac{\pi}{\cos\left[\pi(s-\alpha)/2\right]} \,,
\label{eq:numerator}
\end{eqnarray}
and for the denominator of \eqref{eq:formula} it holds
\begin{eqnarray}
\Gamma\left(\frac{1-s}{2}\right)
\Gamma\left(\frac{1+s}{2}\right)
&=&
\Gamma\left(\frac{1-s}{2}\right)
\Gamma\left(1-\frac{1-s}{2}\right) \nonumber \\
&=& \frac{\pi}{\sin\left[\pi(1-s)/2\right]} \nonumber \\
&=& \frac{\pi}{\cos(\pi s/2)} \,.
\label{eq:denominator}
\end{eqnarray}
Finally, plugging \eqref{eq:numerator} and \eqref{eq:denominator} 
into \eqref{eq:formula} with $n=1$ we have formula \eqref{eq:corollary}. 
Thus, Corollary \ref{Cor1}
is proved.
\proofend

\subsection{Proof of Theorems \ref{Th:DRDC} and \ref{Th:DRDRL} 
and further results}
Let $J^\alpha$ be the Riemann--Liouville fractional integral 
of order $\alpha > 0$ defined by
\be
J^\alpha f(t)= \frac{1}{\Gamma(\alpha)} \int_0^y
(t-\tau)^{\alpha-1} f(\tau) \, d\tau \,,
\quad t \in \RR^+ \,,
\ee
then the Caputo fractional derivatives $\DC_t^\alpha$  
and in the Riemann--Liouville fractional derivative $\DRL_t^\alpha$
of order $m-1 < \alpha < m$, with $m \in \NN$, 
are defined as
\be
\DC_t^\alpha f(t)= J^{m-\alpha} \frac{d^m f}{dt^m} \,,
\label{eq:DC}
\ee
\be
\DRL_t^\alpha f(t) = \frac{d^m}{dt^m} J^{m-\alpha} f(t) \,.
\label{eq:DRL}
\ee
In particular, 
the following noteworthy relation between $\DC_t^\alpha$ and $\DRL_t^\alpha$
exists \cite{mainardi_etal-fcaa-2001}
\be
\DC_t^\alpha f(t)= \DRL_t^\alpha f(t) - 
\sum_{j=0}^{m-1} f^{(j)}(0^+) \frac{t^{j-\alpha}}{\Gamma(j-\alpha+1)} 
\,.
\label{eq:DCDRL}
\ee

Interestingly, 
the Mellin transform of the Riesz fractional derivative 
\eqref{eq:corollary} differs from that of the 
Caputo fractional derivate 
\cite[Section 2.10.4]{podlubny-1999} by a factor built through
trigonometric functions. 
In particular, it holds
\begin{remark}
\label{rem:RC} 
Let $\phi:\RR \to \RR$ be a symmetric function such that
$\phi(x)=\phi(|x|)$ and the following limits are true
\[
i) \, \lim_{|x| \to 0} \phi^{(j)}(|x|) \, |x|^{s-\alpha+j} = 0 \,, \quad
ii) \, \lim_{|x| \to +\infty} \phi^{(j)}|x|^{s-\alpha+j} = 0 \,, 
\]
with $j=0,\dots,m-1$ and proper $\Re(s)$, 
then by multiplying and dividing \eqref{eq:corollary} 
by $\Gamma[1-(s-\alpha)]\Gamma(1-s)$,
by using \eqref{eq:gammaproperty2} and by remembering that
$\sin 2\theta=2\sin\theta\cos\theta$ it holds
\be
\mM(\mD_x^\alpha f)(s) = 
- \frac{\sin[\pi(s-\alpha)/2]}{\sin(\pi s/2)} \,
\mM(\DC_{|x|}^\alpha f)(s) \,, 
\label{eq:RC}
\ee
for $s \ne 0$ and $0 < \Re(s) < 1$.
\end{remark}

Moreover, 
from relation \eqref{eq:DCDRL} we have that 
\begin{remark}
\label{rem:RRL}
The two fractional derivatives in the Caputo \eqref{eq:DC} and
Riemann--Liouville \eqref{eq:DRL} sense 
are equal if $f^{(j)}(0^+)=0$ with $j=0,\dots,m-1$. 
Therefore, from Remark \ref{rem:RC}, within the same conditions,
the following formula holds
\be
\mM(\mD_x^\alpha f)(s) = 
- \frac{\sin[\pi(s-\alpha)/2]}{\sin(\pi s/2)} \,
\mM(\DRL_{|x|}^\alpha f)(s) \,, 
\label{eq:RRL}
\ee
for $s \ne 0$ and $0 < \Re(s) < 1$.
\end{remark}

\begin{remark}
By using the formula $\sin(\theta-\pi)=-\sin \theta$,
from formula \eqref{eq:RC}, as well as from \eqref{eq:RRL}, 
it follows that the two involved Mellin transforms are 
identical when $\alpha=2$. Thus, 
the identity holds also for the two involved operators as it is expected 
for symmetric functions depending on $|x|$.
\end{remark}

By using Remark \ref{rem:RC} we can prove Theorem \ref{Th:DRDC}
as follows.
 
\proof
We observe that formula \eqref{eq:RC} can be understood
as the product of two Mellin transforms. 
In particular, we introduce a function $\mE$ such that
\be
\mM \mE (s)
=-\frac{\sin[\pi(s-\alpha)/2]}{\sin(\pi s/2)} 
=\frac{\sin[\pi(\alpha-s)/2]}{\sin(\pi s/2)} \,,
\label{eq:ME}
\ee
for $0 < \Re(s) < 1$. Thus, 
it results that function $\mE(\xi)$ with $\xi \in \RR^+$ is 
(see the Appendix for calculation)
\be
\mE(\xi)=\frac{2}{\pi}\frac{\sin(\pi\alpha/2)}{1-\xi^2} 
- \cos\left(\frac{\pi\alpha}{2}\right) \, \xi \delta(\xi-1) \,,
\quad \xi \in \RR^+ \,.
\label{eq:E}
\ee
Finally, by using the convolution formula of the Mellin transform,
from \eqref{eq:RC} we have \eqref{eq:DRDC}.
This completes the proof of Theorem \ref{Th:DRDC}.
\proofend

\smallskip
By using Remark \ref{rem:RRL} we can prove Theorem \ref{Th:DRDRL}.

\proof
Following the same reasoning of Theorem \ref{Th:DRDC},
we interpret formula \eqref{eq:RRL} as the product of two
Mellin transforms. Therefore, 
by using again the convolution formula of the Mellin transform,
the same function $\mE$ can be introduced and 
from \eqref{eq:RRL} we have \eqref{eq:DRDRL}. 
This completes the proof of Theorem \ref{Th:DRDRL}.
\proofend

\section{Application to space-fractional diffusion equation}
\label{sec:4}
\setcounter{section}{4} 
\setcounter{equation}{0} 

The space-fractional diffusion equation is a diffusion equation
where the Laplacian is replaced by the fractional Laplacian,
see, e.g., \cite{hanyga-prsla-2001}.
The same holds in the one-dimensional case,
where the second derivative in
space is replaced by the Riesz (or Riesz--Feller) space-fractional
derivative with fractional order $\alpha \in (0,2)$
\cite{gorenflo_etal-fcaa-2000,mainardi_etal-fcaa-2001}.

It is well-known that the solution of the Cauchy problem 
for standard diffusion with initial datum $\delta(x)$
is the Gaussian, or normal, distribution. In this case,
if $t > 0$ is the time variable, the solution is
self-similar with scaling law $t^{-1/2}$,
which is also the law ruling the linear growth in time of the 
standard deviation of particle displacements. 
In the case of space-fractional diffusion,
the solution in the same setting is a L\'evy stable distribution
with scaling law $t^{-1/\alpha}$ with $\alpha \in (0,2)$
and infinite standard deviation.
As a matter of fact, the relation between space-fractional diffusion
and L\'evy stable densities was indeed established by
W. Feller \cite{feller-1952} as reported in Section \ref{sec:1}.

The diffusive processes underlying the space-fractional diffusion
are known as L\'evy flights and they are Markovian processes 
whose jump-size distribution has power-law decaying tails 
\cite{chechkin_etal-anotrans-2008}
and then the walkers' distribution converges to a L\'evy stable density,
see, e.g., 
\cite{shlesinger_etal-jsp-1982,gorenflo_etal-fcaa-1998,valdinoci-bsema-2009}.
Thus, the space-fractional diffusion equation is the governing equation
of L\'evy flights.
These procesess have an important place in physics literature
for modelling anomalous diffusion in disordered media 
\cite{shlesinger_etal-jsp-1982,ott_etal-prl-1990,
bouchaud_etal-pr-1990,metzler_etal-pr-2000,metzler_etal-jpa-2004,
zaburdaev_etal-rmp-2015}. However, 
some issues rise when the experimental implications of
nonlocal models are considered \cite{hilfer-fcaa-2015}.

In spite of the fact that the exact determination of
the jump-size distribution is not necessary for determining the convergence
\cite{shlesinger_etal-jsp-1982,klafter_sokolov-2011,meerschaert2012}, 
the actual timescale for convergence can be indeed 
dependent on the specific form of the jump-size distribution,
as it is the case with some bi-modal distributions 
\cite{pagnini_etal-fcaa-2021}.
This has some implication on models 
in optimal searching or foraging 
that neglect the exact form of the jump distribution,
in particular with respect to transiency and recurrency 
\cite{pagnini_etal-fcaa-2021}.

By assuming a unitary diffusion coefficient with
physical dimension $[L]^2 [T]^{-\alpha}$,
the space-fractional diffusion equation is
\be
\frac{\partial P}{\partial t} = \mL P = -(-\Delta)^{\alpha/2} P \,,
\quad P(x;0)=\delta(x) \,, 
\quad \text{in} \quad \RR^n \times (0,\infty) \,.
\label{eq:SFDE}
\ee
We apply to both sides of \eqref{eq:SFDE} 
the Mellin transform \eqref{eq:mellint}.
From Theorem \ref{Th2}, see \eqref{eq:th2}, 
we have that the LHS reads
\be
\frac{\partial \mM P}{\partial t} =
\pi^{n/2-s}  
\frac{\Gamma\left(\frac{s}{2}\right)} 
{\Gamma\left(\frac{n-s}{2}\right)} 
\frac{\partial}{\partial t} \mM (\mF' f) (n-s) \,,
\label{eq:LHS}
\ee
and from Theorem \ref{Th1}, see \eqref{eq:th1},
we have that the RHS reads
\be
\mM(\mL P)(s)=
- 2^\alpha 
\frac{\Gamma\left(\frac{s}{2}\right) 
\Gamma\left[\frac{n-(s-\alpha)}{2}\right]}
{\Gamma\left(\frac{n-s}{2}\right) 
\Gamma\left[\frac{s-\alpha}{2}\right]}
\mM P(s-\alpha) \,,
\label{eq:RHS}
\ee
and, by applying formula \eqref{eq:th2} of Theorem \ref{Th2}
to the term $\mM P(s-\alpha)$, we finally have
\be
\mM(\mL P)(s)=
- 2^\alpha 
\pi^{n/2-(s-\alpha)}  
\frac{\Gamma\left(\frac{s}{2}\right)}
{\Gamma\left(\frac{n-s}{2}\right)}
\mM (\mF' P) (n-(s-\alpha)) \,.
\label{eq:RHS2}
\ee
By comparing \eqref{eq:LHS} and \eqref{eq:RHS2},
it holds
\be
\frac{\partial}{\partial t} \mM (\mF' P) (n-s) =
- (2\pi)^\alpha 
\mM (\mF' P) (n-(s-\alpha)) \,,
\ee
that can be explicitly written as
\be
\!\!
\int_0^\infty 
\frac{\partial \mF' P}{\partial t} |\kappa|^{(n-s) -1} \, d|\kappa| =
- (2\pi)^\alpha 
\int_0^\infty 
\mF' P(|\kappa|) \, |\kappa|^{(n-(s-\alpha))-1} \, d|\kappa| \,,
\ee 
where on the LHS the Mellin transform and the derivation in time
are exchanged. Finally, by equalling the integrands and solving with
respect to $\mF' P$, it results 
\be
\mF' P(|\kappa|)=e^{-(2\pi|\kappa|)^{\alpha}t} \,,
\ee
that means $P(x;t)=t^{-1/\alpha}\mP(|x|t^{-1/\alpha})$ 
is a radial L\'evy stable distribution 
with stable parameter $\alpha \in (0,2)$.

The same holds for the one-dimensional case,
when the Riesz space-fractional derivative \eqref{RieszD},
i.e., the symmetric Riesz--Feller fractional operator,
is used in the RHS of \eqref{eq:SFDE}
and Corollary \ref{Cor1} through formula \eqref{eq:corollary}
is applied.

\section{Conclusions} 
\label{sec:5}
\setcounter{section}{5} 
\setcounter{equation}{0} 

In this paper we have derived a further definition of the 
fractional Laplacian on the basis of the Mellin transform.
In an unbounded domain, this definition is equivalent to
some other existing definitions whose equivalence has been
already discussed \cite{kwasnicki-fcaa-2017}. 
In general, this Mellin definition can be used for solving space-fractional
differential equation for radial functions.
In particular, we have applied this result in the case of 
the space-fractional diffusion equation, 
that is a remarkable example where the fractional Laplacian is involved.

If the one-dimensional case is considered, 
then the fractional Laplacian reduces to the Riesz 
(i.e., the symmetric Riesz--Feller) fractional derivative and the 
radial property turns into the symmetric property.
In this setting, our results correct a previous 
determination of the Mellin definition of the Riesz pseudo-differential operator
\cite{khan_etal-rp-2021}.
Beside this, we have also established the ralation of the Riesz
fractional derivative with the Caputo and with the Riemann--Liouville
fractional derivatives.

In the view of the recent efforts done for providing 
a rigourous theory of a generalised fractional calculus
\cite{hanyga-fcaa-2020,diethelm_etal-fcaa-2020,
luchko-s-2021,luchko-m-2021,luchko-fcaa-2021,tarasov-m-2023} 
and the corresponding desired future advancements 
\cite{hilfer_etal-m-2019,diethelm_etal-nd-2022},
the contribution of the present research to that program is 
twofold. 

In the case of multi-dimensional operators \cite{tarasov-m-2023},
we underline that settling down the fractional Laplacian,
which is the paradigmatic multi-dimensional nonlocal operator,
in the framework of the Mellin integral transform theory
(Theorem \ref{Th1}) allows for the use of the related Mellin machinery.
This is a strong tool for deriving new results
\cite{marichev-1983,yakubovich_luchko-1994,paris_kaminski-2001}. 
Actually, the success of the Mellin setting in relation with fractional calculus, 
together with the Mellin--Barnes integrals, 
is an established fact both for theoretical results, see,
e.g., \cite{mainardi_etal-fcaa-2003,luchko_etal-fcaa-2013,
bardaro_etal-jfaa-2015}, 
and for resolution techniques of fractional equations, 
see, e.g., \cite{gorenflo_etal-fcaa-2000,mainardi_etal-fcaa-2001}.  

Moreover, when the one-dimensional case is considered,
the Mellin definition of the fractional Laplacian allows for 
the Mellin definition of two-sided fractional derivatives
(Corollary \ref{Cor1}).
This result is new and it joins with well-known results for one-sided derivatives
\cite{podlubny-1999,luchko_etal-fcaa-2013}.
By comparing the Mellin definitions for two-sided and one-sided
fractional derivatives their relationship emerges 
(Theorems \ref{Th:DRDC} and \ref{Th:DRDRL}).
This relationship is a key result that can lead to extend to two-sided operators
the present theory foundation that holds only for one-sided operators
because it is based on the Sonine equations 
\cite{hanyga-fcaa-2020,luchko-fcaa-2021}. 

\section*{Appendix: 
Derivation of the kernel function $\mE(\xi)$ \eqref{eq:E}}

By using the rule
$\sin(\theta+\omega)=\sin\theta\cos\omega + \cos\theta\sin\omega$,
formula \eqref{eq:ME} can be re-written as
\[
\mM\mE(s)=\sin\left(\frac{\pi\alpha}{2}\right)
\, \cot\left(\frac{\pi s}{2}\right) - 
\cos\left(\frac{\pi\alpha}{2}\right) \,,
\]
and then we want two functions $\varphi_1(r)$ and $\varphi_2(r)$
such that $\mM \varphi_1 = \cot (\pi s/2)$ and $\mM \varphi_2 = 1$. 

By noting from \cite[formula (2.24)]{marichev-1983}
that
\[
\int_0^\infty \frac{x^{s-1}}{1-x} \, dx = \pi \cot \pi s \,, 
\quad \Re(s) \in (0,1) \,,
\]
we have, after replacing $x \to r^2$, the first pair
\[
\mM \varphi_1 (s)= \frac{\pi}{2} \cot \frac{\pi s}{2} \,, 
\quad \text{with} \quad \varphi_1(r)=\frac{1}{1-r^2} \,.
\]
Now, let $\mM^{-1}$ be the Mellin inverse transformation
then for a given function $\psi(r)$ it holds $\psi=\mM^{-1}\mM \psi$.
Therefore, after the exchange of integration order we have
\[
\psi(r)=\int_0^\infty \psi(y) \left\{
\mM^{-1} y^s \right\} \frac{dy}{y} \,.
\]
This integral can be understood as a Mellin convolution integral where
by definition $\mM^{-1} y^s = q(r/y)$ and, 
by using the properties of the delta-function, it holds 
$q(z)=z \delta(z-1)$ with $z=r/y$. 
Function $q(r/y)$ is consistent with $\mM q = y^s$.
Thus, by setting $y=1$, we have also the second pair
\[
\mM\varphi_2(s)=1 \,, 
\quad \text{with} \quad \varphi_2(r)=r\delta(r-1) \,.
\]
Finally it results
\[
\mM\mE(s)=
\frac{2}{\pi}\sin\left(\frac{\pi\alpha}{2}\right)
\mM \varphi_1(s) - 
\cos\left(\frac{\pi\alpha}{2}\right) 
\mM \varphi_2(s) \,,
\quad \Re(s) \in (0,1) \,,
\]
from which formula \eqref{eq:E} follows.

However, 
also setting $\mM^{-1} y^s = q^*(r/y)= \delta(r/y-1)$ is appropriate 
for solving the integral though properties of delta-function and 
it is consistent with the Mellin convolution interpretation and
with formula $\mM q^*=y^s$.
Thus, in our case, by setting $y=1$ another suitable choice 
for function $\varphi_2(r)$ is $\varphi_2(r)=\delta(r-1)$.
Moreover, by disregarding the Mellin convolution interpretation,
also setting $\mM^{-1} y^s = q'(r,y)=y \delta(r-y)$ or 
$\mM^{-1} y^s = q{''}(r,y)=r \delta(r-y)$ is appropriate 
and consistent with the formula $\mM q'=\mM q{''}=y^s$.
Actually, this proceedure does not solve the ambiguity emerging from 
different rearrangements of the delta-function 
like $a\delta(r-a)$ and $r\delta(r-a)$, with $a>0$, 
that provide the same Mellin transform.
More extensive analysis on the Mellin transform and delta-function
are available in literature 
\cite{sudland_etal-fcaa-2004,sudland_etal-itsf-2019}.


\begin{acknowledgements}
The authors would like to acknowledge Prof. F. Mainardi and Dr. S. Vitali
for useful discussions and hints that are hidden in the lines of 
the present research. 
GP is supported by the Basque Government through the 
BERC 2022--2025 program and by the Ministry of Science and Innovation: 
BCAM Severo Ochoa accreditation 
CEX2021-001142-S / MICIN / AEI / 10.13039/501100011033.
\end{acknowledgements}

\section*{\small
Conflict of interest} 

{\small
GP is a member of the editorial board of the journal
\textit{Fractional Calculus and Applied Analysis}.
The peer-review process was guided by an independent editor.
The authors declare that they have no other conflict of interest.
}





\bigskip  

\small 
\noindent
{\bf Publisher's Note}
Springer Nature remains neutral with regard to jurisdictional claims in published maps and institutional affiliations.

\end{document}